\documentclass[12pt]{amsart}

\usepackage[USenglish]{babel}
\usepackage{amsmath,amsthm,amssymb,amsfonts}
\usepackage{mathtools}
\usepackage{mathrsfs}
\usepackage{bbm}

\usepackage{indentfirst}
\usepackage{enumerate}
\usepackage{accents}
\usepackage{cancel} 
\usepackage{xifthen}

\usepackage{multirow}
\usepackage{longtable}
\usepackage{enumitem}

\usepackage{float}
\usepackage{graphicx}
\usepackage{placeins}
\usepackage{caption}
\usepackage{arydshln}

\usepackage{tikz}
\usetikzlibrary{patterns}
\usepackage{tikz-cd}

\usepackage{chemgreek,textgreek} 

\usepackage[breaklinks=true, bookmarksopenlevel=1, bookmarksdepth=2]{hyperref}

\allowdisplaybreaks

\newtheorem{thm}{Theorem}[section]

\newtheorem{lem}[thm]{Lemma}

\theoremstyle{plain} 
\newcommand{\thistheoremname}{}
\newtheorem*{genericthm}{\thistheoremname}

\theoremstyle{definition}

\theoremstyle{remark}
\newtheorem{rem}[thm]{Remark}

\numberwithin{equation}{section}

\newcommand{\Z}{\mathbb{Z}}      
\newcommand{\Q}{\mathbb{Q}}      
\newcommand{\R}{\mathbb{R}}      
\newcommand{\C}{\mathbb{C}}      

\newcommand{\eps}{\varepsilon}   

\renewcommand{\L}{\mathcal{L}}
\newcommand{\Lc}{\mathcal{L}^{*}{}}
\newcommand{\selb}{\mathbf{S}}

\newcommand{\pselb}{\mathbf{PS}}
\renewcommand{\k}{m}

\restylefloat{table}

\setlist[itemize]{leftmargin=*}
\setlist[enumerate]{leftmargin=*}

\makeatletter
 \newcommand*\alphgreek[1]{\expandafter\@alphgreek\csname c@#1\endcsname}
 \newcommand*\@alphgreek[1]{\csname chemgreek_int_to_greek:n\endcsname{#1}}
 \newcommand*\Alphgreek[1]{\expandafter\@Alphgreek\csname c@#1\endcsname}
 \newcommand*\@Alphgreek[1]{\csname chemgreek_int_to_Greek:n\endcsname{#1}}

 \AddEnumerateCounter*{\alphgreek}{\@alphgreek}{\chemalpha}
 \AddEnumerateCounter*{\Alphgreek}{\@Alphgreek}{\chemAlpha}
\makeatother

\makeatletter
\def\@tocline#1#2#3#4#5#6#7{\relax
  \ifnum #1>\c@tocdepth 
  \else
    \par \addpenalty\@secpenalty\addvspace{#2}%
    \begingroup \hyphenpenalty\@M
    \@ifempty{#4}{%
      \@tempdima\csname r@tocindent\number#1\endcsname\relax
    }{%
      \@tempdima#4\relax
    }%
    \parindent\z@ \leftskip#3\relax \advance\leftskip\@tempdima\relax
    \rightskip\@pnumwidth plus4em \parfillskip-\@pnumwidth
    #5\leavevmode\hskip-\@tempdima
      \ifcase #1
       \or\or \hskip 1em \or \hskip 2em \else \hskip 3em \fi%
      #6\nobreak\relax
    \dotfill\hbox to\@pnumwidth{\@tocpagenum{#7}}\par
    \nobreak
    \endgroup
  \fi}
\makeatother

\begin{document}

\title[Zeros near \texorpdfstring{$\MakeLowercase{s}=1$}{\MakeLowercase{s}=1} and the constant term of $L'/L$]{Zeros near \texorpdfstring{$\MakeLowercase{s}=1$}{\MakeLowercase{s}=1} and the constant term of $L'/L$ for $L$-functions in the Selberg class}%
\author{Christian T\'afula}%
\address{D\'epartment de Math\'ematiques et Statistique, %
 Universit\'e de Montr\'eal, %
 CP 6128 succ Centre-Ville, %
 Montreal, QC H3C 3J7, Canada}%
\email{christian.tafula.santos@umontreal.ca}%


\begin{abstract}
 Let $\mathcal{L}(s) = \sum_{n=1}^{\infty} a_n n^{-s}$ be an $L$-function in the Selberg class, and $q_{\mathcal{L}}$ its conductor. Let $\ell_0(\mathcal{L})$ be the constant term of the Laurent expansion of $\mathcal{L}'/\mathcal{L}$ at $s=1$. We show that for certain families $\mathcal{F}$ of $L$-functions in the Selberg class with polynomial Euler product:
 \begin{itemize}
  \item If $\mathcal{L}\in\mathcal{F}$ has no zeros $\beta + i\gamma$ with $\beta > 1 - \delta(\log q_{\mathcal{L}})^{-1}$, $|\gamma| < (\log q_{\mathcal{L}})^{-1/2}$ for some absolute $\delta >0$, then $\Re(\ell_0(\mathcal{L})) \ll_{\mathcal{F}} \log q_{\mathcal{L}}$;
  
  \item If $\Re(\ell_0(\mathcal{L})) \ll \log q_{\mathcal{L}}$ for all $\mathcal{L}\in \mathcal{F}$, then there is some absolute $\delta > 0$ such that $\mathcal{L}$ has no zeros $\beta + i\gamma$ with $\beta > 1 - \delta(\log q_{\mathcal{L}})^{-1}$, $|\gamma| < (1-\beta)^{1/2}(\log q_{\mathcal{L}})^{-1/2}$.
 \end{itemize}
 This generalizes, for instance, the case of families of Dedekind zeta functions of number fields with bounded degree.
\end{abstract}
\maketitle
 
\section{Introduction}
 Let $\L(s)$ be an $L$-function in the Selberg class, and 
 \begin{equation}
  \frac{\L'}{\L}(s) = -\frac{\k_{\L}}{s-1} + \ell_0(\L) + \ell_1(\L)(s-1) + \ldots \label{eukro}
 \end{equation}
 be the Laurent expansion of $\L'/\L$ at $s=1$, where $\k_{\L}$ is the order of the pole at $s=1$. For $L$-functions of global fields, $\ell _0(\L)$ is the so-called \emph{Euler--Kronecker constant} of the field (cf. Ihara \cite{ihar06}). Classically, the \emph{non-existence of Siegel zeros} for Dirichlet $L$-functions $L(s,\chi)$, which asks for there not to be a real zero $\beta$ satisfying $1- \delta (\log q)^{-1} < \beta < 1$ for $\chi\pmod{q}$, with $\delta > 0$ fixed and $q$ large, can be shown to be equivalent to the statement ``$|\ell_0(L(\,\cdot\,,\chi))| \ll \log q$''. The same is true for the Dedekind zeta of number fields $\zeta_K$, which becomes ``$|\ell_0(\zeta_K)| \ll_{d_K} \log |\Delta_K|$'', where $\Delta_K$ is the discriminant of $K/\Q$, and $d_K = [K:\Q]$ the degree. This follows from the classical zero-free regions together with estimates for $\frac{\zeta'_K}{\zeta_K}(s) + \frac{1}{s-1}$ in terms of zeros (cf. Lemmas 8.1, 8.2 together with Lemma 5.6 of Lagarias--Odlyzko \cite{lagodl77}).
 
 Our goal in this paper is to extend this relationship for certain natural families $\mathcal{F}$ of $L$-functions in the Selberg class; namely, \emph{admissible families with polynomial Euler product} \eqref{adms} (see Section \ref{defsS} for definitions). Let $q = q_{\L}$ denote the \emph{conductor} of $\L$. For $\L$ in such families $\mathcal{F}$, we prove the following:
 
 \begin{thm}\label{MT1}
  The following hold:
  \begin{enumerate}[label=\textnormal{(\roman*)}]
   \item If $\Re(\ell_0(\L)) \ll \log q_{\L}$ for $\L\in\mathcal{F}$, then there is a constant $\delta\in\R_{>0}$ such that $\L\in\mathcal{F}$ has no zeros in the region
  \[ \left\{s = \sigma + it \in \C ~\Bigg|~ \sigma > 1- \frac{\delta}{\log q_{\L}},\ |t| < \frac{(1-\sigma)^{1/2}}{\sqrt{\log q_{\L}}} \right\}.\medskip \]
  
   \item Let $f:\R_{\geq 1} \to \R$ be a function satisfying $1\leq f(q)\ll \log q$. If $\L\in\mathcal{F}$ has no zeros in the region
  \[ \left\{s = \sigma + it \in \C ~\Bigg|~ \sigma > 1- \frac{1}{f(q_{\L})},\ |t| < \frac{1}{\sqrt{f(q_{\L})}} \right\}, \]
  then $\Re(\ell_0(\L)) \ll_{\mathcal{F}} \sqrt{f(q_{\L})\log q_{\L}}$.
  \end{enumerate}
 \end{thm}
 
 \FloatBarrier
 \begin{figure}[!htb]
 \centering
 \begin{tikzpicture}[scale=6, every node/.style={scale=2}]
   
  \draw[black] (0, -1/2)--(0, 1/2);
  \draw[black, dashed] (0.5, -1/2)--(0.5, 1/2);
  \draw[black] (1, -1/2)--(1, 1/2);
  \draw[black, ->] (-0.5, 0)--(1.05, 0);
      
  \draw[black] (-1/2,0.06/2)--(-1/2,-0.06/2) node[scale=.5, anchor=north] {$-.5\ \ $};
  \draw[black] (0,0.06/2)--(0,-0.06/2) node[scale=.5, anchor=north west] {$0$};
  
  \draw[black] (.5,0.06/2)--(.5,-0.06/2) node[scale=.5, anchor=north west] {$.5$};
  \draw[black,thick] (1-.0005,0.04/2) to[in=87, out=-87] (1-0.0005,-0.04/2) node[scale=.5, anchor=north west] {$1$};
  
  \foreach \xx in {-8,...,15}
   \draw[black] (\xx/16,0.03/2)--(\xx/16,-0.03/2);
  
  \draw[black] (0.03/2, 1/2)--(-0.03/2, 1/2); 
  \draw[black] (0.03/2, -1/2)--(-0.03/2, -1/2); 
  
  \draw[black] (0.9375, 0)--(0.9, .1) node[scale=.5, anchor=south east] {$1-\dfrac{\delta}{\log q_{\L}}$};
  
  \foreach \yy in {-4,...,4}
  \draw[black] (-0.03/2,\yy/8)--(0.03/2,\yy/8);
  
  \draw[scale=1, domain=0.9375:1, smooth, variable=\t, black, thick] plot ({\t},{sqrt(1-\t)});
  \fill[scale=1, domain=0.9375:1, smooth, variable=\t, pattern=north east lines, pattern color=black!50] plot ({\t},{sqrt(1-\t)}) |- (0.9375,0);
  \draw[scale=1, domain=0.9375:1, smooth, variable=\t, black, thick] plot ({\t},{-sqrt(1-\t)});
  \fill[scale=1, domain=0.9375:1, smooth, variable=\t, pattern=north east lines, pattern color=black!50] plot ({\t},{-sqrt(1-\t)}) |- (0.9375,0);
 \end{tikzpicture}
  
  \caption{Shape of the region in Theorem \ref{MT1} (i).}
 \end{figure}
 
 Theorem \ref{MT1} implies that for admissible families $\mathcal{F}$ with polynomial Euler product, $\L\in\mathcal{F}$ has no \emph{Siegel zeros} if $\Re(\ell_0(\L)) \ll_{\mathcal{F}} \log q_{\L}$ -- see Remark \ref{nsiegz}. Examples of such families include Hecke $L$-functions associated to primitive finite-order characters for number fields with bounded degree, and normalized $L$-functions associated with holomorphic modular forms of bounded weight (cf. Kaczorowski--Perelli \cite[p. 317]{kaczpereVI}). For the family $\mathcal{F}$ of real odd Dirichlet characters with very smooth conductor, for instance, the zero-free regions of Theorem \ref{MT1} (ii) hold true for some $f(q) = o_{q\to+\infty}(\log q)$ except for possible Siegel zeros, which do not exist under a certain uniform formulation of the $abc$-conjecture -- cf. \cite[Corollary 1.5]{taf21}.
 

 We prove Theorem \ref{MT1} using a generalization of our previous results in \cite[Section 3]{taf21}, which are Theorems \ref{AnLem} and \ref{AnLem2} below. 
 
\section{Definitions and notation}\label{defsS}
 \subsection{Selberg class}
 A Dirichlet series $\L(s) := \sum_{n\geq 1} a_n n^{-s}$ is said to be an $L$-function in the \emph{Selberg class} $\selb$ if it satisfies the five axioms (cf. Kaczorowski--Perelli \cite{kaczpere99}):
 \begin{enumerate}[label=\textbf{(S\arabic*)}]
  \item\label{s1} $\sum_{n\geq 1} a_n n^{-s}$ is absolutely convergent for $\Re(s)>1$. \smallskip
  
  \item\label{s2} (Analytic continuation) There exists a smallest $\k = \k_{\L}\in\Z_{\geq 0}$ such that $(s-1)^{\k}\L(s)$ is entire of finite order.\smallskip
   
  \item\label{s3} (Functional equation) There exists a triple $(c, Q, \gamma_\L)$, where $c, Q$ are positive real numbers, and $\gamma_\L$ is a function of the form
  \[ \gamma_{\L}(s) := \prod_{j=1}^{f} \Gamma(\lambda_j s + \mu_j) \]
  (called a \emph{Gamma factor}), with $f\in\Z_{\geq 1}$, positive real numbers $\lambda_j$ ($1\leq j\leq f$), and complex numbers $W$, $\mu_j$ ($1\leq j \leq f$) with $|W|=1$ and $\Re(\mu_j)\geq 0$, for which the \emph{completed $L$-function}
  \begin{equation}
   \Lc(s) := c\, Q^s\, \gamma_{\L}(s)\, \L(s) \label{fnceq}
  \end{equation}
  satisfies a functional equation: $\Lc(s) = W \overline{\Lc(1-\overline{s})}$.\smallskip
  
  \item\label{s4} (Ramanujan Hypothesis) $a_n \ll_{\eps} n^{\eps}$ for every fixed $\eps > 0$.\smallskip
  
  \item\label{s5} (Euler product) For $\Re(s) >1$,
  \[ \log \L(s) = \sum_{n\geq 0} \frac{b(n)}{n^{s}}, \]
  where $b(n) =0$ unless $n = p^{k}$, $k\geq 1$ for some positive rational prime $p$, and $b(n) \ll n^{\vartheta}$ for some $0\leq \vartheta <\frac{1}{2}$.
 \end{enumerate}
 
\subsection{Polynomial Euler product}
 We consider a subclass of $\selb$, the \emph{Selberg class with polynomial Euler product} denoted $\pselb$, where axiom \ref{s5} is replaced by axiom \ref{s5p} below:
 \begin{enumerate}[label=\textbf{(S\arabic*{}')}]
  \setcounter{enumi}{4}
  \item\label{s5p} (Polynomial Euler product) There is $r=r_{\L}$ such that, for $\Re(s) >1$ we have
  \[ \L(s) = \prod_{p} \prod_{j=1}^{r}\bigg(1 - \frac{\alpha_j(p)}{p^s}\bigg)^{-1}, \]
  where the product runs over the positive rational primes $p$, and $b(n) = b_{\L}(n) \ll n^{\vartheta}$, where the implied constant is uniform in $\mathcal{F}$.
 \end{enumerate}
 Under \ref{s5p}, the coefficients $a_n$ of the Dirichlet series of $\L$ are multiplicative, and axiom \ref{s4} is equivalent to having $|\alpha_{j}(p)| \leq 1$ ($1\leq j\leq r$) for every $p$ (cf. Steuding \cite[Lemma 2.2]{steuding07}). Thus, since $\log \L(s) = \sum_{p} \sum_{k\geq 1} (\sum_{j=1}^{r}\alpha_{j}(p)^{k})/kp^{ks}$,
 \begin{equation}
  |b(p^k)| \leq \frac{|\alpha_1(p)^k| + \ldots + |\alpha_r(p)^{k}|}{k} \leq \frac{r_{\L}}{k} \label{bddb}
 \end{equation}
 in the notation of \ref{s5}. All automorphic $L$-functions have polynomial Euler product. It is conjectured that $\selb = \pselb$, and that the ``arithmetic'' degree $r_{\L}$ in \ref{s5p} equals the ``analytic'' degree $d_{\L}$ defined below \cite[p. 163]{kacz06}. 
 
 \subsection{Admissible families in \texorpdfstring{$\pselb$}{S'}}
 The $\gamma$-factor in \ref{s3} is not unique; for instance, one can apply the duplication formula to any of the $\Gamma$-functions in the product. Two $\gamma$-factors of a same $\L\in\selb$ differ at most by a constant factor \cite[Theorem 8.1]{kaczpere99}. Quantities that do not change under transformations of the triple $(c,Q,\gamma_{\L})$ are called \emph{invariants} of $\L$. These include
 \begin{equation}
  \log Q + \frac{\gamma'_{\L}}{\gamma_{\L}}(1),\qquad d_{\L} := 2\sum_{j=1}^{f} \lambda_j, \qquad q_{\L} := (2\pi)^{d_{\L}} Q^2 \prod_{j=1}^f \lambda_{j}^{2\lambda_j};
 \end{equation}
 the last two being the \emph{degree} and \emph{conductor}, respectively.
 
 We say $\mathcal{F}\subseteq \pselb$ is an \emph{admissible family} (cf. \cite[p. 317]{kaczpereVI}) if $\L\in\mathcal{F}$ a $\gamma$-factor with
 \begin{equation}
  \begin{gathered}
   Q\gg 1,\quad \lambda_{\min} := \min\limits_{1\leq j\leq f} \lambda_j \gg 1,\quad \max_{1\leq k\leq f} |\mu_{k}| \ll 1; \quad\text{and} \\
   d_{\L} \ll 1,\qquad \k_{\L} \ll 1,\qquad r_{\L} \ll 1,
  \end{gathered}\label{adms}
 \end{equation}
 with $r_{\L}$ as in \ref{s5p}. In such families,
 \begin{equation}
  \bigg|\frac{\gamma'_{\L}}{\gamma_{\L}}(1)\bigg| \ll 1,\qquad \log q_{\L} = 2\log Q + O(1). \label{safql}
 \end{equation}
 so ``$\log Q + \frac{\gamma'_{\L}}{\gamma_{\L}}(1)$'' and ``$\frac{1}{2}\log q_{\L}$'' are essentially interchangeable for large $q_{\L}$.
 
 \begin{rem}
  Conjecturally, it is expected that every $\L\in\selb$ has a $\gamma$-factor with $\lambda_{j} = \frac{1}{2}$ ($1\leq j\leq f$), and that $d_{\L}$, $q_{\L} \in \Z_{\geq 0}$ (cf. \cite[Section 9]{kaczpere99}). 
 \end{rem}
 
\subsection{Zeros of \texorpdfstring{$\L$}{L}}
 Denote by $\varrho(\L)$ the multiset\footnote{A multiset is a set that allows multiple instances of a same element.} of \emph{non-trivial zeros} of $\L$; i.e., zeros of $\Lc(s)$. By the Euler product \ref{s5}, $\L$ never vanishes for $\Re(s)>1$, and by the functional equation \ref{s3}, the only zeros for $\Re(s)$ must correspond to the poles of the $\gamma$-factor, so if $\rho$ is a non-trivial zero then $0\leq \Re(\rho) \leq 1$.
 
 As $\Lc$ has order $1$ as an entire function, $\sum_{\varrho(\L)\neq 0} 1/\varrho$ converges in the principal-value sense ``$\lim\limits_{T\to +\infty} \sum_{\varrho(\L)\neq 0,\, |\Im(\varrho)|\leq T}$''. We have (cf. \cite[Theorem A.1]{taf21}):
 \begin{equation}
  \frac{\L'}{\L}(s) + \frac{\k_{\L}}{s} + \frac{\k_{\L}}{s-1} = \bigg(\sum_{\varrho(\L)} \frac{1}{s-\varrho}\bigg) - \log Q - \frac{\gamma'_{\L}}{\gamma_{\L}}(s). \label{exforF}
 \end{equation}
 
 \begin{rem}[Siegel zeros in the Selberg class]\label{nsiegz}
  Fix a constant $\delta>0$. For any $\L\in\selb$, if $\L(\beta) = 0$ for some real $\beta$ satisfying
  \[ 1 - \frac{\delta}{\log(1 + \lambda Q)} \leq \beta \leq 1, \]
  where $\lambda := \max_{1\leq j\leq f} (\lambda_j + |\mu_j|)$, then $\beta$ is called a \emph{Siegel zero} for $\L(s)$ relative to $\delta$. If $a_n\geq 0$ for every $n\geq 1$, and $\frac{\L'}{\L}(\sigma) \in \R_{<0}$ for every $\sigma >0$, then there is some effectively computable $\delta >0$ for which $\L$ has at most $\k_{\L}$ Siegel zeros (cf. Lemma in Goldfeld--Hoffstein--Lieman \cite{golhoflie94}). From \eqref{safql}, Theorem \ref{MT1} implies that for admissible families $\mathcal{F}\subseteq \pselb$, $\L$ has no Siegel zeros if $\Re(\ell_0(\L)) \ll_{\mathcal{F}} \log q_{\L}$.
 \end{rem}

\section{Lower bounds for \texorpdfstring{$\Re\ell_0(\L)$}{Re l\_0(L)}}
\subsection{Lemmas}
 From \eqref{exforF}, for real $x> 1$ we have
 \begin{equation}
  \sum_{\varrho(\L)} \frac{\varPi_{x-1}(\varrho)}{4} = \Re\left(\frac{\L'}{\L}(x)\right) + \frac{\k_{\L}}{x} + \frac{\k_{\L}}{x-1} + \log Q + \Re\left(\frac{\gamma'_{\L}}{\gamma_{\L}}(x) \right), \label{reLL}
 \end{equation}
 where $\varPi$ is the \emph{pairing function} (cf. \cite[Subsection 3.1]{taf21})
 \begin{equation}
  \varPi_{\eps}(s) := \frac{1}{s + \eps} + \frac{1}{\overline{s} + \eps} + \frac{1}{1-s + \eps} + \frac{1}{1-\overline{s} + \eps}, \label{prupf}
 \end{equation}
 defined for $s,\eps \in \C$ such that $\eps \neq -s, -\overline{s}, -1+s, -1+\overline{s}$. For $0\leq \Re(s) \leq 1$, $\eps \in \R$, one checks that $\varPi_{\eps}(s) \in \R_{> 0}$. Setting $x=1$ in \eqref{reLL}, \eqref{eukro} gives us
 \begin{equation}
  \sum_{\varrho(\L)} \frac{\varPi_{0}(\varrho)}{4} = \Re(\ell_0(\L)) + \k_{\L} + \log Q + \Re\bigg(\frac{\gamma'_{\L}}{\gamma_{\L}}(1)\bigg), \label{pi01}
 \end{equation}
 which is our starting point.
 
 \begin{lem}\label{smllp}
  Let $\L\in\pselb$, and write $\lambda_{\min} = \min_{1\leq j\leq f} \lambda_{j}$. For $\eps >0$ we have:
  \begin{enumerate}[label=\textnormal{(\roman*)}]  
   \item $\displaystyle 0 <\sum_{\varrho(\L)} \frac{\varPi_{\eps}(\varrho)}{4} < \Bigg|\log Q + \Re\left(\frac{\gamma'_{\L}}{\gamma_{\L}}(1) \right)\Bigg| + \frac{2\k_{\L} + r_{\L}}{\eps} + \bigg(\frac{1}{\lambda_{\min}^2} + \frac{\pi^2}{6}\bigg)\frac{d_{\L}}{2}\,\eps$, \smallskip
    
   \item {\small$\displaystyle \Bigg|\sum_{\varrho(\L)} \frac{\varPi_{\eps}(\varrho)}{4} - \k_{\L} - \log Q - \Re\bigg(\frac{\gamma'_{\L}}{\gamma_{\L}}(1)\bigg)\Bigg| < \frac{\k_{\L} + r_{\L}}{\eps} + \bigg(\k_{\L} + \bigg(\frac{1}{\lambda_{\min}^2} + \frac{\pi^2}{6}\bigg)\frac{d_{\L}}{2}\bigg)\eps$}.
  \end{enumerate}
 \end{lem}
 \begin{proof}
  We prove the parts separately.
  
  \medskip
  \noindent
  $\bullet~\text{\underline{Part} (i):}$
  Setting $s = 1+\eps$ in \eqref{reLL} yields\footnote{If one assumes the strong $\lambda$-conjecture (i.e., $\lambda_{j}=\frac{1}{2}$, $1\leq j\leq f$), along with $\mu_j \in \R$, $0\leq \mu_j < .96$ ($1\leq j\leq f$), it is possible to obtain better bounds by having that $\frac{\gamma'_{\L}}{\gamma_{\L}}(1+\eps) < 0$ for small $\eps>0$, which is true for Dirichlet $L$-functions for example (cf. \cite[Lemma 3.3 (i)]{taf21}). This is the case because $\frac{\Gamma'}{\Gamma}(x)$ is strictly increasing for $x\in\R_{>0}$, and $\frac{\Gamma'}{\Gamma}(1.461632\ldots) = 0$.}
  \begin{align}
   \sum_{\varrho(\L)} \frac{\varPi_{\eps}(\varrho)}{4} &\leq \Bigg|\Re\left(\frac{\L'}{\L}(1+\eps)\right)\Bigg| + \frac{2\k_{\L}}{\eps} + \Bigg|\log Q + \Re\left(\frac{\gamma'_{\L}}{\gamma_{\L}}(1+\eps) \right)\Bigg| \label{subst}
  \end{align}  
  Using that $\frac{\Gamma'}{\Gamma}(1+s) = -\gamma + \sum_{k\geq 1} \big(k^{-1} - (k+s)^{-1}\big)$ for $\Re(s) > 0$, where $\gamma = 0.5772\ldots$ is Euler--Mascheroni's constant, we have
  \begin{align*}
   \bigg|\frac{\gamma'_{\L}}{\gamma_{\L}}(1+\eps) &\,- \frac{\gamma'_{\L}}{\gamma_{\L}}(1)\bigg| = \sum_{j=1}^f \lambda_{j}\bigg|\frac{\Gamma'}{\Gamma}(\lambda_j+\mu_j+\eps) - \frac{\Gamma'}{\Gamma}(\lambda_j+\mu_j) \bigg| \\
   &= \eps \sum_{j=1}^f \lambda_j\bigg|\sum_{k= 0}^{\infty} \frac{1}{(k+\lambda_j+\mu_j)(k+\lambda_j+\mu_j+\eps)}\bigg| \\
   &\leq \eps \sum_{j=1}^f \lambda_j\bigg(\frac{1}{\lambda_j(\lambda_j+\eps)} + \sum_{k=1}^{\infty} \frac{1}{n^2}\bigg) = \eps \sum_{j=1}^f \lambda_j\bigg(\frac{1}{\lambda_j(\lambda_j+\eps)} + \frac{\pi^2}{6}\bigg),
  \end{align*}
  and thus
  \begin{equation}
   \bigg|\frac{\gamma'_{\L}}{\gamma_{\L}}(1+\eps) - \frac{\gamma'_{\L}}{\gamma_{\L}}(1)\bigg| < \bigg(\frac{1}{\lambda_{\min}^2} + \frac{\pi^2}{6}\bigg)\frac{d_{\L}}{2}\,\eps. \label{gmmfc}
  \end{equation}
  Moreover, by \eqref{bddb}
  \begin{equation}
   \bigg|\frac{\L'}{\L}(\sigma)\bigg| \leq \sum_{p}\sum_{k\geq 1} \frac{|b(p^k)|\, k\log p}{p^{k\sigma}} \leq r_{\L}\sum_{p}\sum_{k\geq 1} \frac{\log p}{p^{k\sigma}} = r_{\L}\,\bigg|\frac{\zeta'}{\zeta}(\sigma)\bigg| \leq \frac{r_{\L}}{\sigma-1} \label{eq3sth}
  \end{equation}
  for $\sigma>1$. Together with \eqref{gmmfc}, this yields part (i).
  
  \medskip
  \noindent
  $\bullet~\text{\underline{Part} (ii):}$
  From \eqref{reLL}, using \eqref{gmmfc} and \eqref{eq3sth} we get
  \begin{align*}
   &\Bigg|\sum_{\varrho(\L)} \frac{\varPi_{\eps}(\varrho)}{4} - \k_{\L} - \log Q -\Re\bigg(\frac{\gamma'_{\L}}{\gamma_{\L}}(1)\bigg)\Bigg| \\
   &\hspace{10em}< \bigg(\k_{\L} + \bigg(\frac{1}{\lambda_{\min}^2} + \frac{\pi^2}{6}\bigg)\frac{d_{\L}}{2}\bigg)\eps + \frac{\k_{\L}}{\eps} + \Bigg|\Re\bigg(\frac{\L'}{\L}(1+\eps)\bigg)\Bigg| \\
   &\hspace{10em}\leq \bigg(\k_{\L} + \bigg(\frac{1}{\lambda_{\min}^2} + \frac{\pi^2}{6}\bigg)\frac{d_{\L}}{2}\bigg)\eps + \frac{\k_{\L} + r_{\L}}{\eps}. \qedhere
  \end{align*}
 \end{proof}
 
 The following lemma is an adaptation of Lemma 3.4 of \cite{taf21}.
 
 \begin{lem}\label{cuteineq}
  The following hold:
  \begin{enumerate}[label=\textnormal{(\roman*)}]
   \item For $0\leq \Re(s) \leq 1$ and $x \geq \dfrac{1+\sqrt{5}}{2}$, we have $\varPi_{0}(s) \geq \dfrac{\varPi_{x-1}(s)}{2x-1}$; \smallskip
   
   \item Let $M\in\R_{\geq 2}$, and
   \[ \mathcal{B} := \bigg\{s\in \C ~\bigg|~ 1-\frac{1}{M} \leq \Re(s) \leq 1,\ |\Im(t)| \leq \frac{1}{M^{1/2}} \bigg\}. \]
   Then, for $s\in \{0\leq \Re(s)\leq  1\}\setminus (\mathcal{B} \cup (1-\mathcal{B}))$ and $0 < \eps < 1$, we have
   \[ |\varPi_{0}(s) - \varPi_{\eps}(s)| \leq 5M\eps \,\varPi_{\eps}(s).\]
  \end{enumerate}
 \end{lem}
 \begin{proof}
  Part (ii) follows from \cite[Lemma 3.4 (ii)]{taf21}, where strict inequality is proved for $s\in \{0< \Re(s)< 1\}\setminus (\mathcal{B} \cup (1-\mathcal{B}))$. Since $\varPi_{\eps}(s)$ is continuous in $s$ for $s\neq -\eps$, $1+\eps$, the inequality holds in the closure of $\{0< \Re(s)< 1\}\setminus (\mathcal{B} \cup (1-\mathcal{B}))$.
  
  We briefly sketch the argument for part (i). Let $s=\sigma + it$ and $x\geq 1$. Writing $\widetilde{\sigma} := \sigma(1-\sigma)$, we have
  \begin{align}
   \frac{\varPi_{x-1}(s)}{2} &= \frac{x + \sigma -1}{(x+ \sigma-1)^2 + t^2} + \frac{x-\sigma}{(x-\sigma)^2 + t^2} \nonumber \\
   &= \bigg(1 + \frac{\widetilde{\sigma}(1-\widetilde{\sigma}) + 3\widetilde{\sigma} t^2 + x(x-1) \big( 1-2\widetilde{\sigma} - x(x-1) - t^2\big)}{\widetilde{\sigma}^2 + (1 - 2\widetilde{\sigma})t^2 + t^4   +   x(x-1)\big(2\widetilde{\sigma} + x(x-1) + 2t^2\big) }\bigg) \frac{(2x-1)}{1+t^2}. \nonumber
  \end{align}
  Since $x(x-1)(2\widetilde{\sigma} + x(x-1) + 2t^2) > 0$, we get
  \begin{equation*}
   \begin{aligned}
    \text{\small $\displaystyle \varPi_0(s) - \frac{\varPi_{x-1}(s)}{2x-1} \geq 2\bigg(\frac{- 1 + 2\widetilde{\sigma} + x(x-1) + t^2}{\widetilde{\sigma}^2 + (1-2\widetilde{\sigma})t^2 + t^4 + x(x-1)\big(2\widetilde{\sigma} + x(x-1) + 2t^2\big)} \bigg) \frac{x(x-1)}{1+t^2}$}.
   \end{aligned}
  \end{equation*}
  Because $x(x-1) \geq 1$ for $x \geq \frac{1+\sqrt{5}}{2}$, part (i) follows.
 \end{proof}

\subsection{General theorems} 
 Theorem \ref{MT1} will be derived from the following, more general statements.

 \begin{thm}\label{AnLem}
  Let $\L\in\pselb$, and let $S\subseteq \varrho(\L)$ be any finite multiset of non-trivial zeros of $\L$ (including the empty set). Then,
  \begin{align*}
   \Re\Bigg(\sum_{\varrho \in S} \frac{1}{1-\varrho}\Bigg) &< \Re(\ell_0(\L)) + \bigg(1-\frac{1}{\sqrt{5}}\bigg)\,\big|\log Q + \Re\big(\tfrac{\gamma'_{\L}}{\gamma_{\L}}(1)\big)\big| + \bigg(1+\frac{1}{\sqrt{5}}\bigg)\,|S|\,+\\
   &\hspace{.25em}+ \bigg(2 - \frac{1}{\sqrt{5}}\bigg)m_{\L} + \frac{1}{2}\bigg(1 + \frac{1}{\sqrt{5}}\bigg)r_{\L} + \frac{1}{4}\bigg(1 - \frac{1}{\sqrt{5}}\bigg)\bigg(\frac{1}{\lambda_{\min}^2} + \frac{\pi^2}{6}\bigg)d_{\L}.
  \end{align*}
 \end{thm}
 \begin{proof}
  Write $\widetilde{S}$ for the multiset consisting of $\varrho$, $1-\varrho$ for each $\varrho\in S$, so that $|\widetilde{S}| = 2|S|$, and write $H = |\log Q + \Re\big(\frac{\gamma'_{\L}}{\gamma_{\L}}(1)\big)|$. From the definition of the pairing function \eqref{prupf}, we have $\frac{1}{2}\Re((1-\varrho)^{-1}) \leq \varPi_0(\varrho)/4$ and $\varPi_{\eps}(\varrho)/4 \leq \eps^{-1}$, for $0\leq \Re(\varrho) \leq 1$, $\eps>0$, and $\varrho \neq 0,1$. Thus, from \eqref{pi01} and Lemmas \ref{smllp} (ii), \ref{cuteineq} (i), for $x\geq \frac{1+\sqrt{5}}{2}$ we have:
  \begin{align*}
   \Re(&\ell_0(\L)) = \sum_{\varrho\in \widetilde{S}} \frac{\varPi_{0}(\varrho)}{4} + \Bigg(\sum_{\varrho(\L)} \frac{\varPi_{0}(\varrho)}{4} - \k_{\L} - H - \sum_{\varrho\in \widetilde{S}} \frac{\varPi_{0}(\varrho)}{4} \Bigg) \\
   &> \Re\Bigg(\sum_{\varrho \in S} \frac{1}{1-\varrho}\Bigg) + \frac{1}{2x-1}\Bigg(\sum_{\varrho(\L)} \frac{\varPi_{x-1}(\varrho)}{4} - \k_{\L} - H - \sum_{\varrho\in \widetilde{S}} \frac{\varPi_{x-1}(\varrho)}{4}\Bigg)\, - \\
    &\hspace{22.5em} -\, \bigg(1-\frac{1}{2x-1}\bigg)\big(\k_{\L} + H\big) \\
   &> \Re\Bigg(\sum_{\varrho \in S} \frac{1}{1-\varrho}\Bigg) -\frac{1}{2x-1}\bigg(\frac{r_{\L}+\k_{\L}}{x-1} + \bigg(\k_{\L} + \bigg(\frac{1}{\lambda_{\min}^2} + \frac{\pi^2}{6}\bigg)\frac{d_{\L}}{2}\bigg)(x - 1)\bigg) \,- \\
    &\hspace{16em} - \frac{1}{2x-1}\frac{2|S|}{x-1} -\, \bigg(1-\frac{1}{2x-1}\bigg)\big(\k_{\L} + H\big) \\
   &= \Re\Bigg(\sum_{\varrho \in S} \frac{1}{1-\varrho}\Bigg) - \bigg(\frac{x-1}{2x-1}\bigg)\,2\bigg|\log Q + \Re\bigg(\frac{\gamma'_{\L}}{\gamma_{\L}}(1)\bigg)\bigg| - \frac{2|S|}{(2x-1)(x-1)} \,+ \\
   &\hspace{1em} -\, \bigg(\frac{3x^2-6x + 4}{(2x-1)(x-1)}\bigg) \k_{\L} - \frac{1}{(2x-1)(x-1)} r_{\L} + \frac{1}{2}\bigg(\frac{x-1}{2x-1}\bigg)\bigg(\frac{1}{\lambda_{\min}^2} - \frac{\pi^2}{6}\bigg)d_{\L}.
  \end{align*}
  The theorem then follows by taking $x= \frac{1+\sqrt{5}}{2}$.
 \end{proof}

 \begin{thm}\label{AnLem2}
  Let $\L\in\pselb$. Suppose that $|\log Q + \Re(\tfrac{\gamma'_{\L}}{\gamma_{\L}}(1))| \geq 1$, and let $1\leq M \leq C\big|\log Q + \Re\big(\tfrac{\gamma'_{\L}}{\gamma_{\L}}(1)\big)\big|$ for some $C\in\R_{\geq 1}$. Then, in the region
  \begin{equation*}
   \mathcal{B} = \mathcal{B}_{M} := \left\{s \in \C ~\bigg|~ \sigma > 1- \frac{1}{M},\  |t| < \frac{1}{\sqrt{M}} \right\}, 
  \end{equation*}
  we have
  \begin{align*}
   &\Bigg|\Re(\ell_0(\L)) - \Re\Bigg(\sum_{\varrho(\L) \cap \mathcal{B}} \frac{1}{1-\varrho}\Bigg)\Bigg| \\
   &\hspace{4.5em}< 2\,\Big(|\varrho(\L)\cap \mathcal{B}| + O\big((1+C^{1/2})(\k_{\L} + r_{\L})\big)\Big) \sqrt{M \big|\log Q + \Re\big(\tfrac{\gamma'_{\L}}{\gamma_{\L}}(1)\big)\big|} \\
   &\hspace{23em}+ O\Bigg( \frac{(1+\lambda_{\min}^{-2})\,d_{\L}}{\big|\log Q + \Re\big(\tfrac{\gamma'_{\L}}{\gamma_{\L}}(1)\big)\big|} \Bigg).
  \end{align*}
  where the constants implied by $O$ are absolute.
 \end{thm}
 \begin{proof}
  Let $0< \eps < 1$, and write $\widetilde{\mathcal{B}} := \mathcal{B} \cup (1-\mathcal{B})$, from \eqref{pi01} we have
  \begin{align*}
   \Re(\ell_0(\L)) - \Re\Bigg(\sum_{\varrho(\L) \cap \mathcal{B}} \frac{1}{1-\varrho}\Bigg) = &\underbrace{\sum_{\varrho(\L)\setminus \widetilde{\mathcal{B}}}\frac{\varPi_0(\varrho) - \varPi_{\eps}(\varrho)}{4}}_{=:\, S_1} \,- \\
   &-\,\underbrace{\Bigg(\sum_{\varrho(\L)\cap \widetilde{\mathcal{B}}} \frac{\varPi_{\eps}(\varrho)}{4} - \sum_{\varrho(\L) \cap \mathcal{B}} \Re\bigg(\frac{1}{\varrho}\bigg)\Bigg)}_{=:\, S_2} \,+\\
   &+\,\underbrace{\Bigg(\sum_{\varrho(\L)} \frac{\varPi_{\eps}(\varrho)}{4} - \k_{\L} - \log Q -\Re\bigg(\frac{\gamma'_{\L}}{\gamma_{\L}}(1)\bigg)\Bigg)}_{=:\, S_3}.
  \end{align*}
  From Lemma \ref{cuteineq} (ii) we have $|S_1| \leq 5M\eps\, \sum_{\varrho(\L)} \varPi_{\eps}(\varrho)/4$, so from Lemma \ref{smllp} (i) it follows that
  \[ |S_1| < 5M\eps\big|\log Q + \Re\big(\tfrac{\gamma'_{\L}}{\gamma_{\L}}(1)\big)\big| + 5M(2\k_{\L} + r_{\L}) + 5M\cdot Jd_{\L}\, \eps^2, \]
  where $J := \frac{1}{2}(\lambda_{\min}^{-2} + \pi^2/6)$. Next, since $\varPi_{\eps}(\varrho)/4 \leq \eps^{-1}$ for $0\leq \Re(\varrho) \leq 1$, $\varrho\neq 0,1$, and $\Re(1/\varrho) \leq M^{-1} < \eps^{-1}$ for $\varrho \in \mathcal{B}$, we have
  \[ |S_2| < 2 |\varrho(\L)\cap \mathcal{B}|\,\eps^{-1}. \]
  From Lemma \ref{smllp} (ii),
  \[ |S_3| < (\k_{\L} + r_{\L})\eps^{-1} + (\k_{\L} +Jd_{\L})\eps. \]
  Therefore, taking $\eps := 1/\sqrt{M |\log Q + \Re(\tfrac{\gamma'_{\L}}{\gamma_{\L}}(1))|}$ produces
  \begin{align*}
   &\text{\small$\Bigg|\Re(\ell_0(\L)) - \Re\Bigg(\sum_{\varrho(\L) \cap \mathcal{B}} \frac{1}{1-\varrho}\Bigg)\Bigg| < \Bigg(5 + 2|\varrho(\L) \cap \mathcal{B}| + \k_{\L} + r_{\L}\Bigg)\sqrt{M |\log Q + \Re(\tfrac{\gamma'_{\L}}{\gamma_{\L}}(1))|}$} \\
   &\hspace{18.5em}\text{\small$+\ 5M(2\k_{\L} + r_{\L}) + \dfrac{5M\cdot Jd_{\L}}{M |\log Q + \Re(\tfrac{\gamma'_{\L}}{\gamma_{\L}}(1))|}$} \\
   &\hspace{24.5em}\text{\small$+\ \dfrac{\k_{\L}+ Jd_{\L}}{\sqrt{M|\log Q + \Re(\tfrac{\gamma'_{\L}}{\gamma_{\L}}(1))|}}$}
  \end{align*}
  Since $M \leq C|\log Q + \Re(\tfrac{\gamma'_{\L}}{\gamma_{\L}}(1))|$, we have $M \leq C^{1/2}\sqrt{M |\log Q + \Re(\tfrac{\gamma'_{\L}}{\gamma_{\L}}(1))|}$, which concludes the proof.
 \end{proof}
 
 \subsection{Proof of Theorem \ref{MT1}}
  Let $\mathcal{F}\subseteq \pselb$ be an admissible family as in \eqref{adms}.
  
  For part (i), assume that $\Re(\ell_{0}(\L)) \ll \log q_{\L}$ in $\mathcal{F}$. Suppose for the sake of contradiction that for each $j \in \Z_{\geq 1}$, there is $\L_j \in \mathcal{F}$ which has a zero $\rho_j = \beta_j + i\gamma_j$ satisfying $\beta_j > 1 - \delta_j (\log q_{\L})^{-1}$, $|\gamma| < (1-\beta_j)^{1/2}(\log q_{\L})^{-1/2}$ with $\delta_j\to_{j\to\infty} 0$. We have
  \begin{align*}
   \Re\bigg(\frac{1}{1- (\beta_j + i\gamma_j)} \bigg) &= \frac{1}{(1-\beta_j) + \frac{\gamma_j^2}{1-\beta_j}} \geq \frac{\log q_{\L}}{2\delta_j}.
  \end{align*}
  However, by \eqref{safql}, taking $S = \{\rho_j\}$ in Theorem \ref{AnLem} leads to contradiction as $j\to\infty$, and thus no such sequence of $\L_j$s can exist.
  
  For part (ii), taking $M = f(q_{\L})$ in Theorem \ref{AnLem2} yields $\Re(\ell_{0}(\L)) \ll_{\mathcal{F}} \sqrt{f(q_{\L})\log q_{\L}}$, since $\rho(\L)\cap \mathcal{B}_M = \varnothing$ by hypothesis.
  
\bibliographystyle{amsplain}
\bibliography{$HOME/Academie/Recherche/_latex/bibliotheca}%
\end{document}